\documentclass[fleqn]{mat01}
\usepackage{times,mathtimy,amssymb,latexsym}
\begin{document}

\setcounter{page}{307} \firstpage{307}

\def\supp{{\rm Supp\,}}
\def\R{{\mathbb R}}
\def\P{{\mathbb P}}
\def\K{{\mathbb K}}
\def\ptj{\frac{\partial}{\partial x_j}}
\def\pt{\partial}
\def\C{{\mathbb C}}
\def\H{{\mathbb H}}
\def\olmk{{\mathbb M}(\kappa)}
\def\olm{{\bar M}}
\def\lsr{{\Delta_{S(r)}}}
\def\esr{\lambda_1(S(r))}
\def\cpn{{\mathbb C\mathbb P}^n}
 \def\innprod#1#2{\left<#1,#2\right>}
\def\orl{\bar{R}}
\def\og{\bar{\gamma'}}
\def\tv{{\tilde v}}
\def\tw{{\tilde w}}
\def\norm#1{\Vert\hskip.2ex #1 \hskip.2ex\Vert}

\newtheorem{theore}{Theorem}
\renewcommand\thetheore{\arabic{theore}}
\newtheorem{theor}[theore]{\bf Theorem}
\newtheorem{defn}[theore]{\rm DEFINITION}
\newtheorem{remar}[theore]{Remark}


\title{A sharp upper bound for the first eigenvalue of the Laplacian of compact
hypersurfaces in rank-1 symmetric spaces}

\markboth{G Santhanam}{An upper bound for $\lambda_1$}

\author{G SANTHANAM}

\address{Department of Mathematics and Statistics, Indian Institute of
Technology, Kanpur~208~016, India\\
\noindent E-mail: santhana@iitk.ac.in}

\volume{117}

\mon{August}

\parts{3}

\pubyear{2007}

\Date{MS received 10 July 2006; revised 27 September 2006}

\begin{abstract}
Let $M$ be a closed hypersurface in a simply connected rank-1 symmetric space $\olm$.
In this paper, we give an upper bound for the first eigenvalue of the Laplacian of $M$
in terms of the Ricci curvature of $\olm$ and the square of the length of the second
fundamental form of the geodesic spheres with center at the center-of-mass of $M$.
\end{abstract}


\maketitle

\section{Introduction}

Let $(\olmk, \hbox{d} s^2)$ denote the simply connected space form
of constant curvature $\kappa$ where $\kappa=0, 1$ or $-1$ and
dimension $n\geq 2$. Let $M$ be a closed hypersurface of $\olmk$.
When $\kappa=0$ and $M$ a closed hypersurface of $\R^n$,
Bleecker--Weiner \cite{BW} proved that the first eigenvalue
$\lambda_1(M)$ of the Laplace operator of $M$ satisfies the
inequality: $\lambda_1(M)\leq \frac{1}{{\rm vol}(M)}\int_M|A|^2$,
where $|A|^2$ denotes the square of the length of the second
fundamental form of the hypersurface $M$. In \cite{R}, Reilly
improved this inequality to show that
$\lambda_1(M)\leq\frac{n-1}{{\rm vol}(M)}\int_M|H|^2$, where $H$
is the mean curvature of the hypersurface $M$. These inequalities
of Bleecker--Weiner and Reilly are also sharp for geodesic spheres
in $\R^n$. Since then, Reilly's inequality has been extended to
hypersurfaces in other simply connected space forms (see \cite{H}
and \cite{JFG} for details and related results).

While trying to understand these results, we noticed that one can
obtain a similar sharp upper bound for the first eigenvalue
$\lambda_1(M)$, of closed hypersurfaces $M$ in rank-1 symmetric
spaces. Namely $\lambda_{1}(M) \leq \frac{1}{{\rm vol}(M)}\int_{M}
\lambda_{1} (S(r))$ where $\lambda_{1}(S(r))$ is the first
eigenvalue of the geodesic sphere $S(r)$ with center at a point
$p_{0}$ called the center-of-mass of the hypersurface $M$ and
radius $r(q) = d(p_{0}, q)$.

We refer to \cite{B} and \cite{PP} for the basic Riemannian geometry  used in this
paper.

\subsection{\it Statement of results}

To state  our results we need the notion of center-of-mass and result on the existence
of center-of-mass for measurable subsets of $\olm$.

Let $(M, g)$ be a complete Riemannian manifold. For a point $p\in M$, we denote by
$c(p)$, the convexity radius of $(M, g)$ at $p$. For a subset $A\subseteq B(q, c(q))$
for some $q\in M$, we let $CA$ denote the convex hull of $A$. Let $\exp_q\colon
T_qM\to M$ be the exponential map and $X=(x_1, x_2,\dots,x_n)$ be the normal
coordinates centered at $q$. We identify $CA$ with $\exp_q^{-1}(CA)$ and we also
denote $g_q(X, X)$ as $\norm{X}^2_q$ for $X\in T_qM$.

We now state and prove the center-of-mass theorem (see also \cite{BK} and
\cite{ARGS}).

\begin{theor}[\!]\label{thm1}
Let $A$ be a measurable subset of $(M, g)$ contained in $B(q_0,
c(q_0))$ for some point $q_0\in M$. Let $G\colon [0,
2c(q_0)]\to\R$ be a continuous function such that $G$ is positive
on $(0, 2c(q_0))$. Then there exists a point $p\in CA$ such that
\begin{equation*}
\int_A G(\norm{X}_p)X{\rm d}V=0,
\end{equation*}
where $X=(x_1, x_2, \dots, x_n)$ denotes the geodesic normal coordinate system centred
at $p$.
\end{theor}

\begin{proof}
For $q\in CA$, we define
\begin{equation*}
      v(q):=\int_AG(\|{X}\|_q) X {\rm d}V,
\end{equation*}
where $X=(x_1, x_2,\dots, x_n)$ is a geodesic normal coordinate system centred at $q$.

We shall now show that the continuous vector field $v$ points
inward along the boundary $\partial CA$ of $CA$. Then the theorem
follows from the Brouwer's fixed point theorem.

Since $CA$ is convex, it is contained in the half-space
\hbox{$H_q:=\{ X\in T_q M\!\!:g(X, \nu(q))\leq 0\}$} for every
$q\in \partial\/CA$, where $\nu(q)$ denotes the outward normal to
$\partial CA$ at $q$. This implies that $g(v(q), \nu(q))  < 0$ for
all $q\in\partial CA$. Thus $v$ points inward along the
boundary\break of $CA$.

We can find a $\delta > 0$ such that $\exp_q(\delta v(q))\in CA$ for every $q\in CA$.
Then the continuous map $f_v : CA \to CA$ defined by
\begin{equation*}
   f_v(q):= \exp_q(\delta v(q))
\end{equation*}
has a fixed point $p\in CA$ by the Brouwer's fixed point theorem.
Hence $v(p)=0$. This completes the proof of the theorem.\hfill
$\Box$\vspace{-1pc}
\end{proof}

\setcounter{theore}{0}
\begin{defn}$\left.\right.$\vspace{.5pc}

\noindent{\rm The point $p$ in the theorem is called a center-of-mass of the
measurable subset $A$ with respect to the mass distribution function $G$.}
\end{defn}

Before we state our results, we fix some notations that we will be using throughout
this paper.

We let $(\olm, \hbox{d}s^2)$ denote any one of the following rank-1 symmetric spaces:
the round sphere $S^n$ with constant sectional curvature $\frac{1}{4}$, complex
projective space $\C\P^n$, quaternionic projective space $\H\P^n$ and the Cayley
projective plane $Ca\P^2$ with sectional curvature $\frac{1}{4}\leq K_{\olm}\leq 1$ or
their non-compact duals with sectional curvature $-1\leq K_{\olm}\leq -\frac{1}{4}$.
We also write ${\rm dim}\,\olm=kn$, where $k=\dim_\R\K$; $\K=\R$, $\C$, $\H$ or $Ca$.
We also let ${\rm Ric}_{\olm}$ denote the Ricci curvature of $\olm$ and also remark
that ${\rm Ric}_{\olm}$ is constant. The round sphere $S^n$ with constant sectional
curvature $1$ is denoted by $(S^n, Can)$

Given a point $p\in\olm$, we let $S(r)$ denote the geodesic sphere of radius $r$ with
center $p$. Let $\Delta_{S(r)}$ denote the Laplacian of the geodesic sphere $S(r)$
with respect to the induced metric and $\lambda_1(S(r))$ denote the first eigenvalue
of $\Delta_{S(r)}$.

\begin{theor}[\!]\label{thm2}
Let $M$ be a closed hypersurface in a simply connected rank-$1$
symmetric space $(\olm, {\rm d}s^2)$ of compact type. Assume that
$M$ is contained in a ball of radius $\pi/2$. Let $p_0$ be the
center-of-mass of $M$ with respect to the mass distribution
function $G(t)=1/t$. Then
\begin{align}\label{eq1}
\lambda_1(M) &\leq\frac{1}{{\rm vol}(M)}\int_M\lambda_1(S(r))\nonumber\\[.4pc]
&=\frac{1}{{\rm vol}(M)}\int_M(\mid\!A(r)\!\mid^2+{\rm Ric}_{\olm}),
\end{align}
where $r(x):=d(p_0, x)$ is the distance from the point $p_0$ to the point $x$.
Furthermore$,$ equality holds in the above inequality iff $M$ is a geodesic sphere
with center $p_0$.
\end{theor}

\begin{theor}[\!]\label{thm3}
Let $M$ be a closed hypersurface in a simply connected rank-$1$
symmetric space $(\olm, {\rm d}s^2)$ of non-compact type. Let
$p_0$ be the center-of-mass of $M$ with respect to the mass
distribution function $G(t)=1/t$. Then
\begin{align}\label{eq2}
\lambda_1(M) &\leq\frac{1}{{\rm vol}(M)}\int_M\lambda_1(S(r))\nonumber\\[.4pc]
&=\frac{1}{{\rm vol}(M)}\int_M(\mid\!A(r)\!\mid^2+{\rm Ric}_{\olm}),
\end{align}
where $r(x):=d(p_0, x)$ is the distance from the point $p_0$ to the point $x$.
Furthermore$,$ equality holds in the above inequality iff $M$ is a geodesic sphere
with center $p_0$.
\end{theor}

\section{Preliminaries}\label{sec2}

Let $0< r< i(\olm)$. Let $S(r)$ be a geodesic sphere of radius $r$
in $(\olm, \hbox{d}s^2)$ centred at a point $m\in M$. We identify
$S(r)$ with the inverse image $\exp^{-1}(S(r))$ with the metric
$\exp_m^*(\hbox{d}s^2|_{S(r)})$. Let $\lsr$ denote the Laplacian
of $S(r)$ and $\esr$ the first eigenvalue of $\lsr$.

\subsection{\it First eigenvalue of $\lsr$}

In this section we study the first eigenvalue $\esr$ of $\lsr$. This is  also done in
\cite{BB} and \cite{ARGS}.

\subsubsection{\it $(\olm, \hbox{\rm d}s^2)$ of compact type.}\label{2.1.1}
If $k=1$, the geodesic sphere $S(r)$ is  isometric to
$(\olm_{n-1}=S^{n-1}, 4\sin^2(r/2)Can)$. Therefore the first
eigenvalue of $S(r)$ is $\esr=\frac{n-1}{4\sin^2(r/2)}$.

We know that there is a canonical Riemannian submersion \begin{equation*} \Pi: S(r)\to
(\olm_{n-1}, 4\sin^2(r/2) \hbox{d}s^2)
\end{equation*}
with connected totally geodesic fibres, where $\olm_{n-1}$ is the
simply connected compact \hbox{rank-$1$} symmetric space of
dimension $k(n-1)$.

If $k\geq 2$, we can decompose the Laplacian $\lsr$  of $S(r)$ as
\begin{equation}\label{eq3}
         \lsr= \frac{1}{4\cos^2(r/2)}\Delta_{(S^{k-1}, Can)} +
               \frac{1}{4\sin^2(r/2)}\Delta_{(S^{kn-1}, Can)}.
\end{equation}
We also know that the Euclidean co-ordinate functions $X_i$'s, for $1\leq i\leq kn$,
are the first  eigenfunctions of $\Delta_{(S^{kn-1}, Can)}$ corresponding to the first
eigenvalue $kn-1$ (see \cite{BGM} for details). Since the fibres are all totally
geodesic, when we restrict these eigenfunctions to each fibre, they become
eigenfunctions of\vspace{-.3pc}
\begin{equation*}
      \frac{1}{4\cos^2(r/2)}\Delta_{(S^{k-1}, Can)}
\end{equation*}
with eigenvalue
\begin{equation*}
      \frac{k-1}{4\cos^2(r/2)}.
\end{equation*}
Hence we get\vspace{-.3pc}
\begin{equation*}
         \lsr X_i= \left(\frac{kn-1}{4\sin^2(r/2)}+\frac{k-1}{4\cos^2(r/2)}\right)X_i
\end{equation*}
for $1\leq i\leq kn$.

Let $\Delta_{\rm H}$ denote the horizontal Laplacian of the Riemannian submersion.
Since
\begin{equation*}
      \Delta_{\rm H}\mid_{\Pi^*C^{\infty}(\olm_{n-1})}=
            \Pi^*\Delta_{(\olm_{n-1}, \, 4\sin^2(r/2) {\rm d}s^2)}
\end{equation*}
all the eigenfunctions of $\Delta_{(\olm_{(n-1)}, 4\sin^2(r/2) {\rm d}s^2)}$ are also
eigenfunctions of $\lsr$ with the same eigenvalues. In particular, the first non-zero
eigenvalue\vspace{-.3pc}
\begin{equation*}
      \frac{2kn}{4\sin^2(r/2)}
\end{equation*}
occurs as an eigenvalue of $\lsr$ also. Now
\begin{equation*}
   \frac{kn-1}{\sin^2(r/2)}+\frac{k-1}{\cos^2(r/2)}<\frac{2kn}{\sin^2(r/2)}
\end{equation*}
iff\vspace{-.3pc}
\begin{equation*}
      \tan(r/2)<\sqrt{\frac{kn+1}{k-1}}.
\end{equation*}
Since $k\geq 1$ and $n\geq 2$, we see that $\sqrt{\frac{kn+1}{k-1}}>1$. Therefore
$\tan(r/2)<\sqrt{\frac{kn+1}{k-1}}$ for $r<\pi/2$ and consequently the functions
$X_i$, for $1\leq i\leq kn$, are all first eigenfunctions of $\lsr$ with the first
eigenvalue\vspace{-.3pc}
\begin{equation*}
      \esr=\frac{kn-1}{4\sin^2(r/2)} +\frac{k-1}{4\cos^2(r/2)}.
\end{equation*}

\subsubsection{\it $(\olm, {\rm d}s^2)$ of non-compact type.}
We denote by $(\olm)^*$ the compact dual of $\olm$.

If $k=1$, the geodesic sphere $S(r)$ is isometric to
$((\olm_{n-1})^*=S^{n-1}, 4\sinh^2(r/2) Can)$ and hence the first
eigenvalue $\esr$ of $S(r)$ is $\frac{n-1}{4\sinh^2(r/2)}$.

We have the canonical Riemannian submersion\vspace{-.3pc}
\begin{equation}
\Pi :(S(r), \hbox{d}s^2\!\mid_{S(r)})\to ((\olm_{n-1})^*, 4\sinh^2(r/2)
\hbox{d}s^2)\nonumber
\end{equation}
with connected totally geodesic fibres.

If $k\geq 2$, we can decompose the Laplacian $\lsr$ as
\begin{equation*}
\lsr = \frac{-1}{4\cosh^2(r/2)}\Delta_{(S^{k-1}, Can)}
         +\frac{1}{4\sinh^2(r/2)}\Delta_{(S^{kn-1}, Can)}.
\end{equation*}
We know that the euclidean coordinate functions $X_i$'s, for
$1\leq i\leq kn$, are eigenfunctions of $\lsr$ with eigenvalue
\begin{equation*}
\esr=\frac{kn-1}{4\sinh^2(r/2)} - \frac{k-1}{4\cosh^2(r/2)}.
\end{equation*}
Since
\begin{equation*}
\Delta_{\rm H}|_{\Pi^*C^{\infty}((\olm_{n-1})^*)}=
         \Pi^*\Delta_{((\olm_{n-1})^*,  \, 4\sinh^2(r/2){\rm d}s^2)}
\end{equation*}
all the eigenfunctions of $\Delta_{((\olm_{n-1})^*, 4\sinh^2(r/2){\rm d}s^2)}$ are
also eigenfunctions of $\lsr$ with the same eigenvalues. In particular, the first
non-zero eigenvalue
\begin{equation*}
      \frac{2kn}{4\sinh^2(r/2)}
\end{equation*}
occurs as an eigenvalue of $\lsr$ also. Now
\begin{equation*}
     \frac{kn-1}{4\sinh^2(r/2)} - \frac{k-1}{4\cosh^2(r/2)}
\end{equation*}
will be the first non-zero eigenvalue of $\lsr$ so long as
\begin{equation*}
      \frac{kn-1}{4\sinh^2(r/2)} - \frac{k-1}{4\cosh^2(r/2)}
      < \frac{2kn}{4\sinh^2(r/2)}.
\end{equation*}
Since the inequality above is valid for all $r>0$, we get that
\begin{equation*}
      \esr=\frac{kn-1}{4\sinh^2(r/2)} - \frac{k-1}{4\cosh^2(r/2)}
\end{equation*} for all $r>0$.

\subsubsection{\it Geometry of $S(r)$.} We will now relate $\esr$ of $S(r)$ with the square of the length
of the second fundamental form of $S(r)$ and the Ricci curvature
of $\olm$.

Let $\gamma$ be a geodesic starting at the point $p$. Let
$R_{\gamma'(t)}\!\!:T_{\gamma'(t)}M\to T_{\gamma'(t)}M$ be the symmetric endomorphism
defined by $R_{\gamma'(t)}(v):=R(v, \gamma'(t))\gamma'(t)$, where $R$ is the curvature
tensor of $\olm$ .

For $0<t<i(\olm)$, we let $A(t)$ denote the Weingarten map $A(\gamma(t))$ of the
smooth hypersurface $S(t)$ at the point $\gamma(t)$. It is known that these family of
symmetric endomorphisms $A(t)$ satisfy the Riccati equation
\begin{equation*}
A'+A^2+R_{\gamma'}=0
\end{equation*}
along the geodesic $\gamma$. Therefore, by taking the trace of these endomorphisms we
get
\begin{align*}
-{\rm Tr}(A'(r))
&={\rm Tr}(A^2(r))+{\rm Tr}(R_{\gamma'(r)})\\[.3pc]
&=|A(r)|^2+{\rm Ric}_{\olm}(\gamma', \gamma')
\end{align*}
where $|A(r)|^2$ is the square of the length of the second fundamental form of the
geodesic sphere $S(r)$ and ${\rm Ric}_{\olm}$ is the Ricci curvature of $\olm$.

Let $E_2$, $E_3$, $\dots$, $E_{kn}$ be an orthonormal basis of $T_{\gamma(r)}S(r)$
such that the vectors $E_i$, for $2\leq i\leq k$, are tangent to the fibre of the
canonical Riemannian submersion\vspace{-.35pc}
\begin{equation*}
\Pi\!:S(r)\to \begin{cases} (\olm_{n-1}, \hbox{d}s^2)&\text{if $(\olm, \hbox{d}s^2)$
is of compact type}
\\[.18pc]
            (\olm_{n-1}^*, \hbox{d}s^2)&\text{if $(\olm, \hbox{d}s^2)$ is of non-compact
             type}.
\end{cases}
\end{equation*}
Then an easy Jacobi field computation shows that\vspace{-.35pc}
\begin{equation*}
A(r)E_i=\begin{cases}
            \cot r\!E_i
            &\text{if $(\olm, \hbox{d}s^2)$ is of compact type} \\[.15pc]
             \coth r\!E_i &\text{if $(\olm, \hbox{d}s^2)$ is of non-compact type}
            \end{cases}
\end{equation*}
for $2\leq i\leq k$ and\vspace{-.28pc}
\begin{equation*}
A(r)E_i=\begin{cases}
            \frac{1}{2}\cot(r/2) E_i
            &\text{if $(\olm, \hbox{d}s^2)$ is of compact type} \\[.15pc]
             \frac{1}{2}\coth(r/2) E_i &\text{if $(\olm, \hbox{d}s^2)$ is of non-compact type}
            \end{cases}
\end{equation*} for $k+1\leq i\leq kn$. Therefore\vspace{-.1pc}
\begin{align*}
\hskip -4pc \text{Tr}(A(r))=\begin{cases} \frac{k(n-1)}{2}\cot(r/2)+(k-1)\cot
r &\text{if $(\olm, \hbox{d}s^2)$ is of compact type} \\[.15pc]
\frac{k(n-1)}{2}\coth(r/2)+(k-1)\coth r &\text{if $(\olm, \hbox{d}s^2)$ is of
non-compact type}
\end{cases}
\end{align*}
and $-\text{Tr}A'(r)=\lambda_1(S(r))$.\vspace{-.15pc}

\section{Proof of Theorems~\ref{thm2} and \ref{thm3}}

We will now prove Theorems~\ref{thm2} and \ref{thm3}.

Let $M$ be a closed hypersurface in $\olm$ contained in a ball of radius $i(\olm)/2$
where $i(\olm)=\infty$ if $\olm$ is non-compact.

By the center-of-mass theorem~(Theorem \ref{thm1}), there exists a
point $p_0$ such that $\int_M\frac{1}{r}X=0$. Since
$X=(x_1,x_2,\dots,x_{kn})$, we see that $\int_M\frac{1}{r}x_i=0$
for $1\leq i\leq kn$.

Let $f_i:=\frac{x_i}{r}$ for $1\leq i\leq kn$. Since $\int_Mf_i=0$, we can use these
functions $f_i$'s as test functions in the Rayleigh quotient. Therefore\vspace{-.28pc}
\begin{equation*}
\lambda_1(M)\int_Mf_i^2\leq\int_M\norm{\nabla^Mf_i}^2
\end{equation*}
for $1\leq i\leq kn$, where $\nabla^M$ denotes the gradient in $M$. Since
$\sum_{i=1}^{kn}x_i^2=r^2$, we see that $\sum_{i=1}^{kn}f_i^2=1$. Hence\vspace{-.28pc}
\begin{equation*}
\lambda_1(M) \,{\rm vol}(M)\leq\sum_{i=1}^{kn}\int_M\norm{\nabla^Mf_i}^2.
\end{equation*}
If we now show that\vspace{-.28pc}
\begin{equation*}
\sum_{i=1}^{kn}\int_M\norm{\nabla^Mf_i}^2\leq\int_M\esr
\end{equation*}
then we are through.

By Green's identity, $\norm{\nabla^Mf_i}^2=f_i\Delta_Mf_i-\Delta_M(f_i^2)$ for $1\leq
i\leq kn$. Therefore,\vspace{-.25pc}
\begin{equation*}
\int_M\norm{\nabla^Mf_i}^2=\int_Mf_i\Delta_Mf_i-\int_M\Delta_M(f_i^2).
\end{equation*}
Since the boundary of $M$ is empty, it follows from the divergence
theorem that $\int_M\Delta_M(f_i^2)=0$ for $1\leq i\leq kn$. Thus
$\int_M\norm{\nabla^Mf_i}^2=\int_Mf_i\Delta_Mf_i$.

Let us now recall that $\Delta=-\frac{\partial^2}{\partial\nu^2}-
        {\rm Tr}(A)\frac{\partial}{\partial\nu}+\Delta_M$, where $\Delta$ is the
Laplacian of $(\olm, \hbox{d}s^2)$, $A$ is the Weingarten map of the hypersurface $M$
and $\nu$ is the unit outward normal to $M$. Using this identity, we write that
$\Delta_Mf_i=\Delta f_i+\frac{\partial^2}{\partial\nu^2}f_i +{\rm Tr}(A)\langle\nabla
f_i, \nu \rangle$.

We will now compute $\Delta f_i$, $\langle\nabla f_i, \nu\rangle$ and
$\partial^2f_i/\partial\nu^2$.

We decompose the Laplacian $\Delta$ of $\olm$ as\vspace{-.25pc}
\begin{equation*}
\Delta=-\frac{\partial^2}{\partial r^2}-{\rm Tr}(A(r)) \frac{\partial}{\partial
r}+\Delta_{S(r)}
\end{equation*}
along the radial geodesics starting at $p_0$. Now
$\frac{\partial}{\partial r}(\frac{x_i}{r})=0$ for $1\leq i\leq
kn$. Therefore $\Delta f_i=\Delta_{S(r)}f_i$. In \S\ref{sec2}, we
have shown that $\lsr f_i=\esr f_i$ for $0<r<i(\olm)/2$.

Now\vspace{-.25pc}
\begin{align*}
f_i\frac{\partial^2f_i}{\partial\nu^2} &=
    \frac{1}{2}\frac{\partial}{\partial\nu}\langle\nabla(f_i^2), \nu\rangle-
    \langle\nabla f_i, \nu\rangle^2\\[.4pc]
&= \frac{1}{2}\frac{\partial^2}{\partial\nu^2}(f_i^2)-
\left(\frac{\partial\,f_i}{\partial\nu}\right)^2.
\end{align*}
Therefore\vspace{-.25pc}
\begin{equation*}
\sum_{i=1}^{kn}f_i\frac{\partial^2f_i}{\partial\nu^2}
=-\sum_{i=1}^{kn}\left(\frac{\partial\,f_i}{\partial\nu}\right)^2
\end{equation*}
and\vspace{-.25pc}
\begin{align*}
\sum_{i=1}^{kn}\int_Mf_i\frac{\partial^2f_i}{\partial\nu^2}
=-\sum_{i=1}^{kn}\int_M\left(\frac{\partial\,f_i}{\partial\nu}\right)^2.
\end{align*}
Similarly,
\begin{align*}
\sum_{i=1}^{kn}\int_M{\rm Tr}(A)f_i\langle\nabla
    f_i, \nu\rangle & = \frac{1}{2} \int_M{\rm Tr}(A)
    \frac{\partial}{\partial\nu}\left(\sum_{i=1}^{kn}f_i^2\right) \\[.4pc]
& = 0.
\end{align*}

We substitute these quantities in
$\sum_{i=1}^{kn}\int_Mf_i\Delta_M f_i$ to get
\begin{align}
\sum_{i=1}^{kn}\int_Mf_i\Delta_M f_i & =
        \sum_{i=1}^{kn}\int_Mf_i\Delta f_i +
        \int_Mf_i\frac{\partial^2f_i}{\partial\nu^2}
        +\sum_{i=1}^{kn}\int_M{\rm Tr}(A)f_i\langle\nabla
        f_i, \nu\rangle\nonumber\\[.4pc]
&=\int_M\lambda_1(S(r))\sum_{i=1}^{kn}f_i^2-
        \sum_{i=1}^{kn}\int_M\left(\frac{\partial\,f_i}{\partial\nu}\right)^2\nonumber
\end{align}
\begin{align}
\phantom{\sum_{i=1}^{kn}\int_Mf_i\Delta_M f_i}&=\int_M\lambda_1(S(r))-
        \sum_{i=1}^{kn}\int_M\left(\frac{\partial\,f_i}{\partial\nu}\right)^2\nonumber\\[.4pc]
&\leq\int_M\lambda_1(S(r))\label{ineq5}\\[.4pc]
&= \int_M(|A(r)|^2+{\rm Ric}_{\olm})\nonumber.
\end{align}
This completes the proof of the inequality.

For the equality part of the proof, we notice that equality holds
in the above inequality iff $\frac{\partial\,f_i}{\partial\nu}=0$,
for $1\leq i\leq kn$, at all points in $M$. This is true iff the
unit outward normal field is same as the radial vector field
$\nabla r$. Hence the equality holds iff $M$ is a geodesic sphere
of radius $d(p_0, M)$ with center $p_0$. This completes the proof
of theorems~\ref{thm2} and \ref{thm3}.

\setcounter{theore}{0}
\begin{remar}{\rm
The analogue of Bleecker--Weiner \cite{BW} result is not known in
rank-1 symmetric spaces. However we have analogous  results of
Theorems \ref{thm2} and \ref{thm3}  in $\R^n$.}
\end{remar}

\setcounter{theore}{3}
\begin{theor}[\!]\label{thm4}
Let $M$ be a closed hypersurface in $\R^n$. Let $p_0$ be the center-of-mass of $M$
with respect to the mass distribution function $G(t)=1/t$. Then
\begin{align}\label{eq6}
\lambda_1(M) &\leq\frac{1}{{\rm vol}(M)}\int_M\lambda_1(S(r))\nonumber\\[.4pc]
&=\frac{1}{{\rm vol}(M)}\int_M|A(r)|^2\nonumber\\[.4pc]
&=\frac{n-1}{{\rm vol}(M)}\int_M\frac{1}{r^2(x)},
\end{align}
where $r(x):=d(p_0, x)$ is the distance from the point $p_0$ to the point $x$.
Furthermore$,$ equality holds in the above inequality iff $M$ is a geodesic sphere
with center $p_0$.
\end{theor}

\begin{proof}
Let us first observe that the inequality~\ref{ineq5} in the proof of
theorems~\ref{thm2} and \ref{thm3} is valid in all Riemannian manifolds in which the
functions $x_i/r$ are first eigenfunctions of the geodesic sphere $S(r)$. This is true
in $\R^n$. Further, the  first eigenvalue $\esr$ of the geodesic sphere $S(r)$ in
$\R^n$ is $\frac{n-1}{r^2}$. Therefore
\begin{align*}
\sum_{i=1}^{kn}\int_Mf_i\Delta_M f_i
&\leq\int_M\lambda_1(S(r))\label{eq5}\\[.4pc]
&=(n-1) \int_M \frac{1}{r^2(x)}.
\end{align*}
Hence
\begin{equation*}
\lambda_1(M)\leq \frac{n-1}{{\rm vol}(M)}\int_M\frac{1}{r^2(x)}.
\end{equation*}

$\left.\right.$\vspace{-1.8pc}

\hfill $\Box$
\end{proof}

\setcounter{theore}{1}
\begin{remar}{\rm
The analogue of Bleecker--Weiner \cite{BW} result in rank-1
symmetric spaces and their relation with the results of
Theorems~\ref{thm2}, \ref{thm3} and \ref{thm4} will be discussed
in a subsequent paper.}
\end{remar}


\begin{thebibliography}{99}
\bibitem{ARGS} Aithal~A~R and Santhanam~G, Sharp upper bound for the
first non-zero Neumann eigenvalue for bounded domains in rank-1 symmetric spaces, {\it
Trans. Am. Math. Soc.} {\bf 348(10)} (1996) 3955--3965

\bibitem{BW} Bleecker~D and Weiner~J, Extrinsic bounds on $\lambda_1$ of
$\Delta$ on a compact manifold, {\it Comment. Math. Helv.} {\bf 51} (1976) 601--609

\bibitem{BB} Bourguignon~J P and Berard Bergery~L, Laplacians
and Riemannian submersions with totally geodesic fibres, {\it
Illinois J. Math.} {\bf 26} (1982) 181--200

\bibitem{B} Berger~M, Lectures on geodesics in Riemannian
geometry (Mumbai: Tata Institute of Fundamental Research) (1965)

\bibitem{BGM} Berger~M, Gauduchon~P and Mazet~E, Le Spectre d'une vari\'{e}t\'{e}
riemannienne, Lecture Notes in Mathematics 194 (Berlin: Springer) (1971)

\bibitem{BK} Buser~P and Karcher~H, Gromov's almost flat
manifolds (Soci\'{e}t\'{e} Math\'{e}matique de France) (1981)

\bibitem{H} Ernst Heintze, Extinsic upper bounds for
$\lambda_1$, {\it Math. Ann.} {\bf 280} (1988) 389--402

\bibitem{JFG} Grosjean~J~F, Upper bounds for the first eigenvalue of the
Laplacian on compact submanifolds, {\it Pacific. J. Math.} {\bf 206} (2002) 93--112

\bibitem{PP} Petersen~P, Riemannian Geometry, GTM Series 171 (New~York:
Springer) (1998)

\bibitem{R} Reilly~R, On the first eigenvalue of the Laplacian for compact submanifold
of Euclidean Space, {\it Comment. Math. Helv.} {\bf 52} (1977) 525--533
\end{thebibliography}
\end{document}